\documentstyle{article}
\def\a{\alpha}    \def\b{\beta}      \def\g{\gamma}
\def\d{\delta}    \def\l{\lambda}    \def\m{\mu}

\def\osp{\widehat{{\rm osp}}(1|2)}

\def\ZZ{Z\!\!\!\!\!\!Z}
\def\CC{{\bf C}\!\!\!{\rm l}}
\newcommand{\bn}{\begin{equation}}
\newcommand{\ed}{\end{equation}}

\newtheorem{proposition}{Proposition}[section]

\normalbaselines
\def\beli{\begin{list}{(\roman{miski})}{\usecounter{miski}}}
\def\eli{\end{list}}
\begin{document}

\begin{center}
\vspace*{1.0cm}

\def\thefootnote{\fnsymbol{footnote}}

{\LARGE{\bf Cartan-Weyl basis for
Quantum Affine Superalgebra $U_{q}(\osp)$ $^*$}}\\
\footnotetext[1]{To be published in the Proceedings of VI-th 
Colloquium
,,Quantum Groups and Integrable Systems'', June 19-21, 1997, Prague, 
Czech
Republic.}
\def\thefootnote{\arabic{footnote}}
\setcounter{footnote}{1}
\vskip 0.6cm
{\large {\bf J. Lukierski$^{1,3}$ and V.N. Tolstoy$^{1,2,5}$}}

\vskip 0.3cm
{\it $^{1}$Institute of Theoretical Physics, University of Wroclaw\\
50-204 Wroclaw \& Poland (e-mail: lukier@ift.uni. wroc.pl)}\\

\vskip 0.1cm
{\it $^{2}$Institute of Nuclear Physics, Moscow State University \\
119899 Moscow \& Russia (e-mail: tolstoy@anna19.npi.msu.su)}
{$^{3}$}
\footnotetext[3]{permanent address.}
\footnotetext[4]{Supported by KBN grant No 2P 30208706.}
\footnotetext[5]{Supported by Russian Foundation for Fundamental 
Research,
grant No 96-01-01421.}

\end{center}

\begin{abstract}
Cartan-Weyl basis for
the quantum affine superalgebra $U_{q}(\osp)$ is constructed in an
explicit form. 
\end{abstract}
\setcounter{equation}{0}
\section{Prelinimaries}
{\large {\it 1. Affine Lie superalgebra $B(0,1)^{(1)}$}}.
The affine untwisted Lie superalgebra (or in other words "Kac - Moody 
affine
Lie superalgebra") $B(0,1)^{(1)}$ has the rank 2, and its Dynkin
diagram is presented by the picture \cite{K2}

\vskip 15pt
\centerline{
\begin{picture}(430,20)
\put(193,9){\line(1,0){32}}
\put(194,7){\line(1,0){32}}
\put(190,6){\circle{8}}
\put(230,6){\circle*{8}}
\put(223,13){\line(1,-2){5}}
\put(223,-1){\line(1,2){5}}
\put(194,5){\line(1,0){32}}
\put(193,3){\line(1,0){32}}
\put(186,14){$\alpha_{0}$}
\put(227,14){$\alpha$}
\end{picture}}
\vskip 8pt
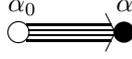
\centerline{\footnotesize Fig.1. Dynkin diagram of the Lie
superalgebra $B(0,1)^{(1)}$}

\vskip 10pt
\noindent
where the dark root $\a$ is odd and the white root $\a_{0}$
is even. The corresponding standard $A=(a_{ij})$ Cartan matrix has
the form:
\bn
A=\pmatrix{{\ }\,\,2&-1\cr-4&{\ }\,\, 2}~.
\label{P7}
\ed
The determinant of this matrix is equal to zero: $\det A=0$.
Moreover we shall use the symmetric $A^{sym}=(a_{ij}^{sym})_{i,j=0,1}$
Cartan matrix:
\bn
A^{sym}=\pmatrix{(\a_{0},\a_{0}) & (\a_{0},\a)\cr
\;(\a,\a_{0})&\;(\a,\a)}=
\pmatrix{{\ }\;4(\a,\a) & -2(\a,\a)\cr
-2(\a,\a) & {\ }{\ }\;\;(\a,\a)}\,,
\label{P8}
\ed
and a non-singular extended symmetric Cartan matrix
$\bar{A}^{sym}=(\bar{a}^{sym}_{ij})_{i,j=-1}^{1}$ and also
its inverse $(\bar{A}^{sym})^{-1}=\big(d_{ij}\big)_{i,j=-1}^{1}$:
\bn
\bar{A}^{sym}=\pmatrix{0 & \;\;\;1 & \;\;\;0 \cr
1 & {\ }\;4(\a,\a) & -2(\a,\a)\cr
0 & -2(\a,\a) & {\ }{\ }\;\;(\a,\a)}\,, \qquad\quad
\big(\bar{A}^{sym}\big)^{-1}=\pmatrix{0 & 1 & 2 \cr
1 & 0 & 0 \cr
2 & 0 & \frac{1}{(\a,\a)}}~.
\label{P9}
\ed
In the case of affine Lie (super)algebras it is convenient to
number rows and columns of the extended symmetric Cartan matrix
$\bar{A}^{sym}=(\bar{a}^{sym}_{ij})_{i,j=-1}^{r}$ by
$-1,0,\ldots,r$ in such way that the first row and column are
numbered beginning with -1 then 0 and so on. Here we keep this
convention. It is also useful to note that the matrix elements of
the standard and symmetric matrices are connected by the relation
\bn
a_{ij}=\frac{2a^{sym}_{ij}}{a^{sym}_{ii}}~.
\label{P10}
\ed

\noindent
The Lie superalgebra $B(0,1)^{(1)}$ is generated by the Chevalley
basis \{$h_{\rm d}$, $h_{\a}$, $h_{\a_{0}}$, $e_{\pm\a}$, 
$e_{\pm\a_{0}}$\}
with the defining relations
\begin{eqnarray}
&&
[h_{\g},h_{\g'}]=0 \qquad\qquad\qquad\;\;
(\g,\;\g'={\rm d},\;\a_{0},\;\a)\,,
\label{P11}
\\
\cr
&&
[e_{\b},e_{-\b'}]=\d_{\b,\b'}h_{\b}\qquad\qquad
(\b,\;\b'=\a,\;\a_{0})\,,
\label{P12}
\\
\cr
&&
[h_{\rm d},e_{\pm\a_{0}}]=\pm e_{\pm\a_{0}}\,,\qquad\quad
[h_{\rm d},e_{\pm\a}]=0\,,
\label{P13}
\end{eqnarray}
\begin{eqnarray}
&&[h_{\a_{0}},e_{\pm\a_{0}}]=\pm 4(\a,\a)e_{\pm\a_{0}}\,,\qquad\;
[h_{\a_{0}},e_{\pm\a}]=\mp 2(\a,\a)e_{\pm\a}\,,
\label{P14}
\\
\cr
&&
[h_{\a},e_{\pm\a_{0}}]=\mp2(\a,\a)e_{\pm\a_{0}}\,,\qquad\;\;\,
[h_{\a},e_{\pm\a}]=\pm (\a,\a)e_{\pm\a}\,,
\label{P15}
\end{eqnarray}
\begin{eqnarray}
&&[e_{\pm\a},[e_{\pm\a},[e_{\pm\a},[e_{\pm\a},
[e_{\pm\a},e_{\pm\a_{0}}]]]]]=0\,,
\label{P16}
\\
\cr
&&[[e_{\pm\a},e_{\pm\a_{0}}],e_{\pm\a_{0}}]=0~.
\label{P17}
\end{eqnarray}
Here the supercommutator $[\cdot,\cdot]$ is given by the formula
\bn
[a,b]=ab-(-1)^{\deg a \deg b } b a\,,
\ed
 where
\bn
\deg h_{\rm d}=\deg h_{\a}=\deg h_{\a_{0}}=\deg e_{\pm\a_{0}}=0, 
\qquad
\deg e_{\pm\a}=1~.
\label{P18}
\ed
It should be noted that a subalgebra generated by the elements
$h_{\a}, e_{\pm\a}$ is isomorphic to $osp(1|2)$.
A subalgebra  of $B(0,1)^{(1)}$ generated by the elements
$h_{\rm d}, h_{\a_{0}}, h_{\a}$ is called the Cartan subalgebra
${\cal H}$ of affine Lie superalgebra $B(0,1)^{(1)}$.
It is easy to see that the element
\bn
h_{\d}:=h_{\a_{0}}+2h_{\a}
\label{P19}
\ed
is a center of superalgebra $B(0,1)^{(1)}$, i.e. it commutes
with every element of $B(0,1)^{(1)}$
\bn
[h_{\d},{\rm everything}]=0~.
\label{P20}
\ed
From (\ref{P19}) we have that
\bn
h_{\a_{0}}=h_{\d}-2h_{\a}=h_{\d-2\a}~.
\label{P21}
\ed
Moreover taking (\ref{P13})-(\ref{P15}) into account we also
obtain that
\bn
(\d,\a)=(\d,\d)=({\rm d},\a)=0\,,\qquad\qquad
(\d,{\rm d})=1~.
\label{P22}
\ed
Frequently we shall also use the notation $\d-2\a$ instead of
$\a_{0}$. The root $\d=\a_{0}+2\a$ is called
the minimal positive imaginary root of $B(0,1)^{(1)}$.

{\large {\it 2. Current realization of $B(0,1)^{(1)}$}}.
One of remarkable properties of affine Lie algebras is that they
admit an exact current or loop realization, i.e. there is an
isomorphism between affine Lie algebras and some extension of
current algebras over finite-dimensional simple Lie algebras.
In particular this property allows easily to obtain the total list
of roots and commutation relations of Cartan - Weyl basis for
the affine Lie algebras. This property is also valid for the supercase
$B(0,1)^{(1)}$. We briefly remind the current realization
(see \cite{K2}).

Let ${\cal L}(g):=g[u,u^{-1}]$ be a linear space of Laurent
polynomials in variable $u$ with coefficients from a simple
finite-dimensional contragredient Lie (super)algebra $g$,
i.e. any element $x\in {\cal L}(g)$ has the form
$x=\sum_{n\in\,\ZZ}\,a_{n}u^{n}$ where all but a finite number of
$a_{k}\in g$ are 0. We specify a Lie (super)algebra structure on
${\cal L}(g)$ with the (super)bracket $[\cdot,\cdot]_{0}$ defined as
follows
\bn
[au^{n},bu^{m}]_{0}=
[a,b]\,u^{n+m} \qquad (a,b\in g;\;\; n,m \in\ZZ\;)~.
\label{P30}
\ed
\noindent
It should be noted that the $\;\ZZ\,_{2}$-grading function on
${\cal L}(g)$ is defined by the condition
\bn
\deg au^{n}=\deg a~.
\label{P31}
\ed
\noindent

Let $(\cdot\mid\cdot)$ be a non-degenerate invariant supersymmetric
bilinear $\CC\,$-valued form on $g$. In the case of superalgebra
$osp(1|2)$ such form is the Killing form and it is unique up to
a constant multiple. We extend this form by linearity to an
${\bf C}[u,u^{-1}]$-valued bilinear form $(\cdot|\cdot)_{\cal L}$
on ${\cal L}(g)$ by
\bn
(au^{n}|bu^{m})_{\cal L}=(a|b)u^{n+m}
\qquad (a,b\in g;\;\;n,m \in \ZZ\;)~.
\label{P32}
\ed
We also define the derivation $\frac{{\rm d}}{{{\rm d}}u}$ of the Lie
(super)algebra ${\cal L}(g)$ by
\bn
\frac{{\rm d}}{{\rm d}u}(au^{n})=a(u^{n})'=nau^{n-1}
\qquad (a\in g;\;\;n\in \ZZ\;)~.
\label{P33}
\ed

Now one can define a $\CC\,$-valued 2-cocycle $\psi (\cdot,\cdot)$
on the Lie (super)algebra ${\cal L}(g)$ as follows:
\bn
\psi (x,y)=
{\rm Res}\left(\frac{\rm d}{{\rm d}u}(x)|y\right)_{\cal L}~.
\label{P34}
\ed
where the linear function (residue) ${\rm Res}$ is
defined by usual way, i.e. ${\rm Res}\;p(u)=c_{-1}$ if
$p(u)=\sum_{k}c_{k}u^{k}$ $(c_{k}\in \CC\,)$. In particular,
for $x=au^{n}$, $b=bu^{m}$ $(a,b \in g;\; n,m\in \ZZ\,)$ we have
\bn
\psi (au^{n},bu^{m}):=n\d_{n,-m}(a|b)~.
\label{P35}
\ed

Recall that  $\CC\,$-valued  2-cocycle on a Lie (super)algebra
${\cal A}$ is a bilinear $\CC\,$-valued function $\psi(\cdot,\cdot)$
satisfying two conditions:
\bn
\psi(x,y)=(-1)^{1+\deg x\deg y}\psi(y,x)\,,
\label{P36}
\ed
\bn
(-1)^{\deg x \deg z}\psi([x,y],z)+(-1)^{\deg z\deg x}\psi([z,y],x)
+(-1)^{\deg x\deg y}\psi([x,z],y)=0~.
\label{P37}
\ed
for any homogeneous elements $x,y,z \in {\cal A}$.
It is not difficult to check that our function (\ref{P34})
satisfies these conditions.

Denote by $\tilde{\cal L}(g)$ the extension of the Lie (super)algebra
${\cal L}(g)$ by a 1-dimensional center, associated to the cocycle
$\psi$. Explicitly the (super)algebra Lie $\tilde{\cal L}(g)$ is
the direct sum of vector spaces:
\bn
\tilde{\cal L}(g)={\cal L}(g)\oplus\CC\,c~.
\label{P38}
\ed
with the (super)bracket defined by
\bn
[au^n+\l c, bu^m+\m c]=[a,b]u^{n+m}+n\d_{n,-m}\psi(a|b)c
\quad(a,b\in g;\;n,m \in\ZZ\,; \;\l,\m \in \CC\;)~.
\label{P39}
\ed
Finally we denote by $\hat{\cal L}(g)$ the Lie (super)algebra which
is obtain by adjoint to $\tilde{\cal L}(g)$ a derivation $d$ which 
acts on
${\cal L}(g)$ as $u\frac{{\rm d}}{{\rm d}u}$ and kills $c$.
More explicitly, $\hat{\cal L}(g)$ is a complex vector space
\bn
\hat{\cal L}(g)=\tilde{\cal L}(g)\oplus\CC\,c\oplus\CC\,d~.
\label{P40}
\ed
with the (super)brackets defined by
\bn
[au^{n}+\l_{1}c+\m_{1}d, bu^{m}+\l_{2}c+\m_{2}d]=
[a,b]u^{n+m}+n\d_{n,-m}(a|b)c+\m_{1}mbu^{m}-\m_{2}nau^{n}.
\label{P41}
\ed
for any
$a,b \in g;\;\;n,m \in \ZZ\,;\;\l_{1},\m_{1},\l_{2},\m_{2} \in 
\CC\;)$~.
The (super)algebra $\hat{\cal L}(g)$ is often denoted as $\hat{g}$
and if it is isomorphic to $g^{(1)}$ then $\hat{g}$ is called
the current realization of $g^{(1)}$.

In our case, i.e. when $g=osp(1|2)$, we can show that the superalgebra
$B(0,1)^{(1)}$ is isomorphic to $\osp$. The isomorphism is
arranged as
\begin{eqnarray}
&&h_{\rm d}\mapsto d\,,\qquad\qquad\;\;
h_{\d}\mapsto c\,,\qquad\qquad\;
h_{\a}\mapsto h_{\a}\,,
\label{P42}
\\
\cr
&&e_{\pm\a}\mapsto e_{\pm\a}\,,\qquad\quad
e_{\pm(\d-2\a)}\mapsto\mp\sqrt{\frac{2}{(\a,\a)}}e_{-2\a}u^{\pm 1}~.
\label{P43}
\end{eqnarray}
This isomorphism immediately allows us to describe the total root 
system and the root space decomposition of $B(0,1)^{(1)}$ 
with respect to the Cartan subalgebra ${\cal H}$ since such 
description is easy obtained for $\osp$.
Indeed, the elements $h_{\a}u^{n}$, $e_{\pm\a}u^{n}$,
$e_{\pm2\a}u^{n}$ for all $n\in \ZZ\,$ are weight
vectors with respect to ${\cal H}=\CC\,h_a\oplus \CC\,c\oplus \CC\,d$
and they form a basis in ${\cal L}(osp(1|2))$. Let $\d$ be a linear
function on ${\cal H}$ defined by $\d\mid_{\CC\,h_{\a}+\CC\,c}=0$,
$\d(d)=1$ then the total root system of $\osp$ has the form
\bn
\Delta:=\left\{n\d\pm\a;\,n\d\pm 2\a;\,n\d\,(n \neq 0)\mid
n\in\ZZ\,\right\}
\label{P44}
\ed
and its root space decomposition can be presented as
\bn
\osp={\cal H}\oplus\left(\bigoplus_{\g}\osp_{\g}\right)\,,
\label{P45}
\ed
where
\begin{eqnarray}
&&\osp_{n\d\pm\a}=\CC\,e_{\pm\a}u^{n}\,,\qquad
\osp_{n\d\pm 2\a}=\CC\,e_{\pm2\a}u^{n}~ \;\;\;\;(n\in\ZZ\,)\,,
\label{P46}
\\
\cr
&&\osp_{n\d}=\CC\,h_{\a}u^{n}\qquad\qquad(n\in\ZZ\,\setminus\{0\})~.
\label{P47}
\end{eqnarray}
We set
\begin{eqnarray}
&&e_{n\d\pm 2\a}=(-1)^{n}\sqrt{\frac{2}{(\a,\a)}}e_{\pm2\a}u^n\,,
\qquad\;\;
e_{-n\d\mp 2\a}=\sqrt{\frac{2}{(\a,\a)}}e_{\mp 2\a}u^{-n}\,,
\label{P48}
\\
\cr
&&e_{n\d\pm\a}=\pm (-1)^{n}e_{\pm\a}u^n\,,\qquad\qquad\qquad
e_{-n\d\mp\a}=\pm e_{-\a}u^{-n}\,,
\label{P49}
\\
\cr
&&e_{n\d}=\frac{(-1)^{n+1}}{\sqrt{(\a,\a)}} h_{\a}u^n\,,
\qquad\qquad\qquad\quad\,
e_{-n\d}=\frac{-1}{\sqrt{(\a,\a)}}h_{\a}u^{-n}~.
\label{P50}
\end{eqnarray}
where $n\in\ZZ\,_{+}\setminus\{0\}$.
The formulas (\ref{P48})-(\ref{P50}) are called
the current realization of the affine root vectors of
$B(0,1)^{(1)}$. The root vectors (\ref{P48})-(\ref{P50})
satisfy the relations (\ref{CW3})-(\ref{CW15}) where $q=1$.
In particular we have
\bn
\begin{array}{rclrcl}
[e_{n\d-2\a},e_{-n\d+2\a}]&=&(-1)^{n+1}h_{n\d-2\a}, &\quad
[e_{n\d+2\a},e_{-n\d-2\a}]&=&(-1)^{n+1}h_{n\d+2\a},
\\
\cr
{}[e_{n\d-\a},e_{-n\d+\a}]&=&(-1)^{n+1}h_{n\d-\a}, &
[e_{n\d+\a},e_{-n\d-\a}]&=&(-1)^{n}h_{n\d+\a},
\\
\cr
[e_{n\d},e_{-n\d}]&=&(-1)^{n}h_{n\d}\,,
\end{array}
\label{P53}
\ed
where $n\in\ZZ\,_{+}\setminus\{0\}$.

Now we return to the root system (\ref{P42}).
The roots $n\d$ $(n\in\ZZ\,\setminus\{0\}$) are called imaginary
roots, the remaining ones are called the real roots.
Moreover we have $\Delta=\Delta_0+\Delta_1$, where
$\Delta_0=\{\pm2\a;\,n\d\pm2\a,\,n\d\mid n\in\ZZ\,\setminus\{0\}\}$
is the even root system and $\Delta_1=\{n\d\pm\a\mid n\in \ZZ\,\}$ is
the odd root system. The positive root system is
\bn
\Delta_{+}:=
\left\{\a,\,2\a;\,n\d\pm\a,\,n\d\pm2\a,\,n\d
\mid n\in\ZZ\,_{+}\setminus\{0\}\right\}~.
\label{P54}
\ed
We denote by $\underline\Delta_+$ the reduced system of positive
roots, which is obtained from $\Delta_+$ by removing the double real
roots $2\a$, $2n\d\pm2\a$ $(n\in\ZZ\,_+\setminus\{0\})$.
It is convenient to present the total root system $\Delta$ by
the picture \cite{L}:
\vskip 10pt
\centerline{
\unitlength 0.29 mm
\begin{picture}(430,120)
\put(16,107){\footnotesize$-3\d+2\alpha$}
\put(0,100){\ldots}
\put(40,100){\circle{5}}
\put(68,107){\footnotesize$-2\d+2\alpha$}
\put(90,100){\circle{5}}
\put(121,107){\footnotesize$-\d+2\alpha$}
\put(140,100){\circle{5}}
\put(188,107){\footnotesize$2\alpha$}
\put(190,100){\circle{5}}
\put(226,107){\footnotesize$\d+2\alpha$}
\put(240,100){\circle{5}}
\put(275,107){\footnotesize$2\d+2\alpha$}
\put(290,100){\circle{5}}
\put(325,107){\footnotesize$3\d+2\alpha$}
\put(340,100){\circle{5}}
\put(375,107){\footnotesize$4\d+2\alpha$}
\put(390,100){\circle{5}}
\put(420,100){\ldots}
\put(16,81){\footnotesize$-3\d+\alpha$}
\put(0,75){\ldots}
\put(40,75){\circle*{5}}
\put(68,81){\footnotesize$-2\d+\alpha$}
\put(90,75){\circle*{5}}
\put(121,81){\footnotesize$-\d+\alpha$}
\put(140,75){\circle*{5}}
\put(188,81){\footnotesize$\alpha$}
\put(190,75){\circle*{5}}
\put(226,81){\footnotesize$\d+\alpha$}
\put(240,75){\circle*{5}}
\put(275,81){\footnotesize$2\d+\alpha$}
\put(290,75){\circle*{5}}
\put(325,81){\footnotesize$3\d+\alpha$}
\put(340,75){\circle*{5}}
\put(375,81){\footnotesize$4\d+\alpha$}
\put(390,75){\circle*{5}}
\put(420,75){\ldots}
\put(33,56){\footnotesize$-3\d$}
\put(0,50){\ldots}
\put(40,50){\circle{5}}
\put(83,56){\footnotesize$-2\d$}
\put(90,50){\circle{5}}
\put(133,56){\footnotesize$-\d$}
\put(140,50){\circle{5}}
\put(190,50){\vector(0,1){22}}
\put(190,50){\vector(1,-1){48}}
\put(237,56){\footnotesize$\d$}
\put(240,50){\circle{5}}
\put(285,56){\footnotesize$2\d$}
\put(290,50){\circle{5}}
\put(335,56){\footnotesize$3\d$}
\put(340,50){\circle{5}}
\put(385,56){\footnotesize$4\d$}
\put(390,50){\circle{5}}
\put(420,50){\ldots}
\put(16,30){\footnotesize$-3\d-\alpha$}
\put(0,25){\ldots}
\put(40,25){\circle*{5}}
\put(69,30){\footnotesize$-2\d-\alpha$}
\put(90,25){\circle*{5}}
\put(120,30){\footnotesize$-\d-\alpha$}
\put(140,25){\circle*{5}}
\put(183,30){\footnotesize$-\alpha$}
\put(190,25){\circle*{5}}
\put(226,30){\footnotesize$\d-\alpha$}
\put(240,25){\circle*{5}}
\put(275,30){\footnotesize$2\d-\alpha$}
\put(290,25){\circle*{5}}
\put(325,30){\footnotesize$3\d-\alpha$}
\put(340,25){\circle*{5}}
\put(375,30){\footnotesize$4\d-\alpha$}
\put(390,25){\circle*{5}}
\put(420,25){\ldots}
\put(16,5){\footnotesize$-3\d-2\alpha$}
\put(0,0){\ldots}
\put(40,0){\circle{5}}
\put(69,5){\footnotesize$-2\d-2\alpha$}
\put(90,0){\circle{5}}
\put(120,5){\footnotesize$-\d-2\alpha$}
\put(140,0){\circle{5}}
\put(182,5){\footnotesize$-2\alpha$}
\put(190,0){\circle{5}}
\put(237,5){\footnotesize$\d-2\alpha$}
\put(240,0){\circle{5}}
\put(275,5){\footnotesize$2\d-2\alpha$}
\put(290,0){\circle{5}}
\put(325,5){\footnotesize$3\d-2\alpha$}
\put(340,0){\circle{5}}
\put(375,5){\footnotesize$4\d-2\alpha$}
\put(390,0){\circle{5}}
\put(420,0){\ldots}
\end{picture}}

\vskip 15pt
\centerline{\footnotesize Fig.2. The total root system
$\Delta$ of $\osp$}
\vskip 12pt
The picture of the reduced root system looks as follows:
\vskip 10pt
\centerline{
\unitlength 0.29 mm
\begin{picture}(430,120)
\put(16,107){\footnotesize$-3\d+2\alpha$}
\put(0,100){\ldots}
\put(40,100){\circle{5}}
\put(121,107){\footnotesize$-\d+2\alpha$}
\put(140,100){\circle{5}}
\put(226,107){\footnotesize$\d+2\alpha$}
\put(240,100){\circle{5}}
\put(325,107){\footnotesize$3\d+2\alpha$}
\put(340,100){\circle{5}}
\put(420,100){\ldots}
\put(16,81){\footnotesize$-3\d+\alpha$}
\put(0,75){\ldots}
\put(40,75){\circle*{5}}
\put(68,81){\footnotesize$-2\d+\alpha$}
\put(90,75){\circle*{5}}
\put(121,81){\footnotesize$-\d+\alpha$}
\put(140,75){\circle*{5}}
\put(188,81){\footnotesize$\alpha$}
\put(190,75){\circle*{5}}
\put(226,81){\footnotesize$\d+\alpha$}
\put(240,75){\circle*{5}}
\put(275,81){\footnotesize$2\d+\alpha$}
\put(290,75){\circle*{5}}
\put(325,81){\footnotesize$3\d+\alpha$}
\put(340,75){\circle*{5}}
\put(375,81){\footnotesize$4\d+\alpha$}
\put(390,75){\circle*{5}}
\put(420,75){\ldots}
\put(33,56){\footnotesize$-3\d$}
\put(0,50){\ldots}
\put(40,50){\circle{5}}
\put(83,56){\footnotesize$-2\d$}
\put(90,50){\circle{5}}
\put(133,56){\footnotesize$-\d$}
\put(140,50){\circle{5}}
\put(190,50){\vector(0,1){22}}
\put(190,50){\vector(1,-1){48}}
\put(237,56){\footnotesize$\d$}
\put(240,50){\circle{5}}
\put(285,56){\footnotesize$2\d$}
\put(290,50){\circle{5}}
\put(335,56){\footnotesize$3\d$}
\put(340,50){\circle{5}}
\put(385,56){\footnotesize$4\d$}
\put(390,50){\circle{5}}
\put(420,50){\ldots}
\put(16,30){\footnotesize$-3\d-\alpha$}
\put(0,25){\ldots}
\put(40,25){\circle*{5}}
\put(69,30){\footnotesize$-2\d-\alpha$}
\put(90,25){\circle*{5}}
\put(120,30){\footnotesize$-\d-\alpha$}
\put(140,25){\circle*{5}}
\put(183,30){\footnotesize$-\alpha$}
\put(190,25){\circle*{5}}
\put(226,30){\footnotesize$\d-\alpha$}
\put(240,25){\circle*{5}}
\put(275,30){\footnotesize$2\d-\alpha$}
\put(290,25){\circle*{5}}
\put(325,30){\footnotesize$3\d-\alpha$}
\put(340,25){\circle*{5}}
\put(375,30){\footnotesize$4\d-\alpha$}
\put(390,25){\circle*{5}}
\put(420,25){\ldots}
\put(16,5){\footnotesize$-3\d-2\alpha$}
\put(0,0){\ldots}
\put(40,0){\circle{5}}
\put(120,5){\footnotesize$-\d-2\alpha$}
\put(140,0){\circle{5}}
\put(237,5){\footnotesize$\d-2\alpha$}
\put(240,0){\circle{5}}
\put(325,5){\footnotesize$3\d-2\alpha$}
\put(340,0){\circle{5}}
\put(420,0){\ldots}
\end{picture}}

\vskip 15pt
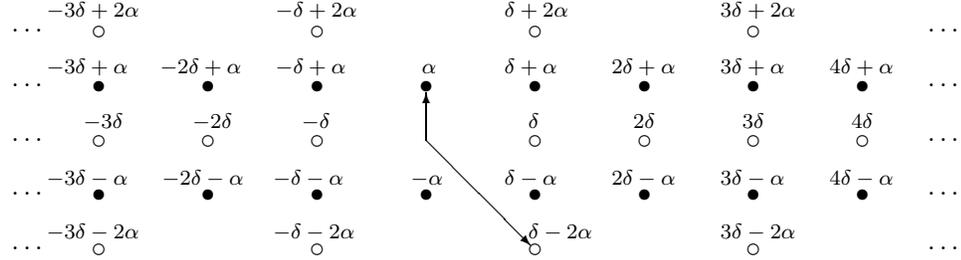
\centerline{\footnotesize Fig.3. The total reduced root system
$\underline{\Delta}$ of $\osp$}

\vskip 10pt

\setcounter{equation}{0}
\section{Defining relations of $U_q(\osp)$}

The $q$-deformation of the universal enveloping algebra $U(\osp)$
(or in other words the quantum affine superalgebra $U_q(\osp)$) is
generated by the Chevalley elements
$k_{\rm d}^{\pm1}:=q^{\pm h_{\rm d}}$,
$k_{\a}^{\pm1}:=q^{\pm h_{\a}}$,
$k_{\d-2\a}^{\pm1}:=q^{\pm h_{\d-2\a}}$,
$e_{\pm\a}$, $e_{\pm(\d-2a)}$
with the defining relations \cite{T1,KT1}
\begin{eqnarray}
&&
k_{\g}k_{\g}^{-1}=k_{\g}^{-1}k_{\g}=1\qquad\qquad\quad\;\;
(\g={\rm d},\,\a,\,\d-2\a)\,,
\label{DR1}
\\
\cr
&&
[k_{\g}^{\pm}, k_{\g'}^{\pm}]=0\qquad\qquad\qquad\qquad\;
(\g,\,\g'={\rm d},\,\a,\,\d-2\a)\,,
\label{DR2}
\\
\cr
&&
k_{\g}e_{\pm\b}k_{\g}^{-1}=q^{\pm(\g,\b)}e_{\pm \b}\qquad\qquad\,
(\b=\a,\,\d-2\a;\;\g={\rm d},\,\b)\,,
\label{DR3}
\\
\cr
&&
[e_{\b},e_{-\b'}]=
\d_{\b,\b'}\frac{k_{\b}-k_{\b}^{-1}}{q-q^{-1}}\qquad\quad
(\b,\;\b'=\a,\;\d-2\a)\,,
\label{DR4}
\\
\cr
&&
[e_{\pm\a},[e_{\pm\a},[e_{\pm\a},[e_{\pm\a},
[e_{\pm\a},e_{\pm(\d-2a)}]_{q}]_{q}]_{q}]_{q}]_{q}=0\,,
\label{DR5}
\\
\cr
&&
[[e_{\pm\a},e_{\pm(\d-2a)}]_{q},e_{\pm(\d-2a)}]_{q}=0\,,
\label{DR6}
\end{eqnarray}
where the $q$-bracket $[\cdot,\cdot]_{q}$ means
the following $q$-supercommutator
\bn
[e_{\b},e_{\b'}]_{q}=e_{\b}e_{\b'}-
(-1)^{\vartheta(\b)\vartheta(\b')} q^{(\b,\b')}e_{\b'}e_{\b}\,,
\label{DR7}
\ed
where the following notation is used
\bn
\vartheta(\b):=\deg e_{\b}~.
\label{DR8}
\ed

A Hopf structure of the quantum superalgebra $U_q(\osp)$ is given
by the following formulas for a comultiplication $\Delta_q$,
an antipode $S_q$ and a counit $\epsilon_{q}$:
\begin{eqnarray}
&&\Delta_{q}(k_\g^{\pm1})=k_\g^{\pm1}\otimes k_\g^{\pm1}\,,
\qquad
\Delta_{q}(e_{\b})=
e_{\b}\otimes 1+k_{\b}^{-1}\otimes e_{\b}\,,
\label{DR9}
\\
\cr
&&\Delta_{q}(e_{-\b})=
e_{-\b}\otimes k_{\b}+1 \otimes e_{-\b}\,,
\label{DR10}
\\
\cr
&&S_{q}(k_\g^{\pm1})=k_\g^{\mp1}\,,\qquad
S_{q}(e_{\b})=-k_{\b}e_{\b}\,,\qquad
S_{q}(e_{-\b})=-e_{-\b}k_{\b}^{-1}\,,
\label{DR11}
\\
\cr
&&\epsilon_{q}(k_\g^{\pm1})=1\,,\qquad
\epsilon_{q}(e_{\pm\b})=0\,,
\label{DR12}
\end{eqnarray}
where $\b=\a,\,\d-2\a$; $\g={\rm d},\,\a,\,\d-2\a$.

\setcounter{equation}{0}
\section{Cartan-Weyl basis for $U_q(\osp)$}

A general scheme for construction of  Cartan-Weyl basis for
quantized Lie algebras and superalgebras was proposed in \cite{T1}.
The scheme was applied in detail at first for quantized
finite-dimensional Lie (super)algebras \cite{KT1} and then to
quantized non-twisted affine algebras \cite{TK}.

This procedure is bases on a notion of ``normal ordering'' for the
reduced positive root system. For affine Lie (super)algebras this
notation was formulated in \cite{T2} (see also \cite{T1,KT2,KT3}).
In our case the reduced positive system has only two normal orderings:
{\small
\bn
\begin{array}{l}
\a,\d+2\a,\d+\a,3\d+2\a,2\d+\a,\ldots,\infty\d+\a,
(2\infty+1)\d+2\a,(\infty+1)\d+\a,\d,2\d,\ldots,\\
\\
\infty\d,(\infty+1)n\d-\a,(2\infty+1)\d-2\a,\infty\d-\a,\ldots,
2\d-\a,3\d-2\a,\d-\a,\d-2\a,
\end{array}
\label{CW1}
\ed}

{\small
\bn
\begin{array}{l}
\d-2\a,\d-\a,3\d-2\a,2\d-\a,\ldots,\infty\d-\a,(2\infty+1)\d-2\a,
(\infty+1)n\d-\a,\d,2\d,\ldots,\\ \\
\infty,(\infty+1)\d+\a,(2\infty+1)\d+2\a,\infty\d+\a,\ldots,
2\d+\a,3\d+2\a,\d+\a,\d+2\a,\a.
\end{array}
\label{CW2}
\ed}

\noindent
The first normal ordering (\ref{CW1} corresponds to ``clockwise''
ordering of positive roots on Fig.3 starting from root $\a$
to root $\d-2\a$. The inverse normal ordering (\ref{CW2})
corresponds to ``anticlockwise'' ordering for the positive roots
when we move from $\d-2\a$ to $\a$.

We choose the normal ordering (\ref{CW1} and in accordance with the
procedure \cite {T1,KT1,TK} we put
\begin{eqnarray}
&&e_{\d-\a}:=\frac{1}{\sqrt{a}}[e_\a,e_{\d-2\a}]_q\,,
\qquad\qquad\quad\;\,
e_{-\d+\a}:=\frac{1}{\sqrt{a}}[e_{-d+2\a},e_{-\a}]_{q^{-1}}\,,
\label{CW3}
\\
\cr
&&e_{\d}:=\frac{1}{\sqrt{b}}[e_{\a},e_{\d-2\a}]_q\,,
\qquad\qquad\qquad\;\,
e_{-\d}:=\frac{1}{\sqrt{b}}[e_{-\d+\a},e_{-\a}]_{q^{-1}}\,,
\label{CW4}
\\
\cr
&&e_{n\d+\a}:=\frac{1}{\sqrt{b}}[e_{(n-1)\d+\a},e_{\d}]\,,
\qquad\qquad
e_{-n\d-\a}:=\frac{1}{\sqrt{b}}[e_{-\d},e_{-(n-1)\d-\a}]\,,
\label{CW5}
\\
\cr
&&e_{(n+1)\d-\a}:=\frac{1}{\sqrt{b}}[e_{\d},e_{n\d-\a}]\,,
\qquad\qquad
e_{-(n+1)\d+\a}:=\frac{1}{\sqrt{b}}[e_{-n\d+\a},e_{-\d}]\,,
\label{CW6}
\\
\cr
&&e_{n\d}'=\frac{1}{\sqrt{b}}[e_{\a},e_{n\d-\a}]_q\,,
\qquad\qquad\qquad\;
e_{-n\d}'=\frac{1}{\sqrt{b}}[e_{-n\d+\a},e_{-\a}]_{q^{-1}}\,,
\label{CW7}
\end{eqnarray}
\begin{eqnarray}
&&e_{(2n-1)\d+2\a}:=
\frac{1}{\sqrt{a}}[e_{(n-1)\d+\a},e_ {n\d+\a}]_{q}\,,
\label{CW8}
\\
\cr
&&e_{-(2n-1)\d-2\a}:=
\frac{1}{\sqrt{a}}[e_{-n\d-\a},e_{-(n-1)\d-\a}]_{q^{-1}}\,,
\label{CW9}
\\
\cr
&& e_{(2n+1)\d-2\a}:=
\frac{1}{\sqrt{a}}[e_{(n+1)\d-\a},e_{n\d-\a}]_{q}\,,
\label{CW10}
\\
\cr
&&e_{-(2n+1)\d+2\a}:=
\frac{1}{\sqrt{a}}[e_{-n\d+\a},e_{-(n+1)\d+\a}]_{q^{-1}}\,,
\label{CW11}
\end{eqnarray}
for $n=1,\,2,\ldots$, where $a=[2(\alpha,\alpha)]$ and 
$b=[2(\alpha,\alpha)]-[(\alpha,\alpha)]$.

The constructed set of the root vectors (\ref{CW3})-(\ref{CW9})
(together with the Chevalley basis) is called a q-analog Cartan-Weyl
basis of $U(\osp)$ or a Cartan-Weyl basis of the quantum algebra
$U_{q}(\osp)$. Now we present the commutation relations
for the vectors of this basis. 
\begin{proposition}
The root vectors $e_{\pm(n\d+\a)}$, $e_{\pm(n\d-\a)}$,
$e_{\pm((2n-1)\d-\a)}$, and $e_{\pm((2n+1)\d+2\a)}$
satisfy the following commutation relations
\begin{eqnarray}
&&[e_{n\d+\a},e_{-n\d-\a}]=
(-1)^{n}\frac{k_{n\d+\a}-k_{n\d+\a}^{-1}}{q-q^{-1}}\,,
\label{CW12}
\\
\cr
&&[e_{n\d-\a},e_{-n\d+\a}]=
(-1)^{n-1}\frac{k_{n\d-\a}-k_{n\d-\a}^{-1}}{q-q^{-1}}\,,
\label{CW13}
\\
\cr
&&[e_{(2n-1)\d+2\a},e_{-(2n-1)\d-2\a}]=
\frac{k_{(2n-1)\d+2\a}-k_{(2n-1)\d+2\a}^{-1}}{q-q^{-1}}\,,
\label{CW14}
\\
\cr
&&[e_{(2n+1)\d-2\a},e_{-(2n+1)\d+2\a}]=
\frac{k_{(2n+1)\d-2\a}-k_{(2n+1)\d-2\a}^{-1}}{q-q^{-1}}\,,
\label{CW15}
\end{eqnarray}
where $n=1,\,2,\ldots$.
\label{CWP1}
\end{proposition}
The sector with imaginary root vectors requires the redefinitions
of the generators (\ref{CW7}). 
We introduce the new imaginary roots vectors $e_{\pm n\d}$
by the following (Schur) relations (see \cite{KT2}, \cite{KT3}
\cite{TK}):
\bn
{e'}_{n\d}=\sum_{p_{1}+2p_{2}+\ldots +np_{n}=n}
\frac{{(q-q^{-1})}^{\sum p_{i}-1}}{p_{1}!\cdots
p_{n}!}{(e_{\d})}^{p_1}\cdots {(e_{n\d})}^{p_n}.
\label{CW16}
\ed
In terms of generating functions
\begin{eqnarray}
&&{\cal E}'(u):=(q-q^{-1})\sum_{n\geq 1}^{}{e'}_{n\d}u^{-n}\,,
\label{CW17}
\\
\cr
&&{\cal E}(u)=(q-q^{-1})\sum_{n\geq 1}^{}e_{n\d}u^{-n}
\label{CW18}
\end{eqnarray}
the relation (\ref{CW16}) may be rewritten in the form
\bn
{\cal E}'(u)=-1+\exp{\cal E}(u)
\label{CW19}
\ed
or
\bn
{\cal E}(u)=\ln(1+{\cal E}'(u))~.
\label{CW20}
\ed
This provides the formula inverse to (\ref{CW16})
\bn
e_{n\d}=\sum_{p_{1}+2p_{2}+\ldots +np_{n}=n}
\frac{{(q^{-1}-q)}^{\sum p_{i}-1} (\sum_{i=1}^{n}p_{i}-1)!}
{p_{1}!\cdots p_{n}!}{{e'}_{\d}^{p_1}
\cdots{e'}_{n\d}^{p_n}}~.
\label{CW21}
\ed
The root vectors of negative roots are obtained by the Cartan
conjugation $(^*)$:
\bn
e_{-n\d}=(e_{n\d})^{*}~.
\label{CW22}
\ed
\begin{proposition}
The root vectors $e_{\pm n\d}$ satisfy the following commutation
relation
\bn
[e_{n\d},e_{-m\d}]=\d_{nm}
\frac{(-1)^{n}}{nb}\left(q^{n(\a,\a)}+q^{n(\a,\a)}-1\right)
[n(\a,\a)]\frac{k_{n\d}-k_{n\d}^{-1}}{q-q^{-1}}\,,
\label{CW23}
\ed
where $n,\,m=1,\,2,\ldots$.
\label{PCW2}
\end{proposition}
{\it Remark}. The root vectors $e_{\pm n\d}'$ do not satisfy the
relation (\ref{CW18}), i.e. $[e_{n\d},e_{-m\d}]\ne 0$ for $n\ne m$.

\setcounter{equation}{0}
\section{Final Remarks}
In this note we restricted ourselves to the description of Cartan-
Weyl basis.
Related results ($q$-Weyl group, automorphisms introducing real forms,
universal $R$-matrix, 
classical $r$-matrix) will be published in our subsequent publication.

\end{document}